\documentstyle[amssymb]{amsart}

\newcommand{\Proof}{\noindent {\sc Proof} \hspace{0.2in}} 

\newtheorem{thm}{Theorem} 
\newtheorem{claim}{Claim}[section]
\newtheorem{lemma}[claim]{Lemma}

\newtheorem{definition}[claim]{Definition}

\theoremstyle{remark}
\newtheorem{remark}{Remark}

\newcommand{\bbr}{{\Bbb R}}

\newcommand{\ood}{{}^{\textstyle \omega\!>} 2}

\newcommand{\omod}{{{}^{\textstyle \omega_{1}}2}}
\newcommand{\omood}{{{}^{\textstyle \omega_{1}\!>}2}}

\newcommand{\mmmd}{{{}^{\textstyle m}2}}

\newcommand{\kkd}{{{}^{\textstyle k}2}}

\newcommand{\od}{{{}^{\textstyle \omega}2}}

\newcommand{\con}{2^{\aleph_0}}

\newcommand{\rng}{\operatorname{rang}}
\newcommand{\cf}{\operatorname{cf}}
\newcommand{\lh}{\operatorname{lh}}
\newcommand{\dom}{\operatorname{dom}}

\newcommand{\otp}{\operatorname{otp}}

\newcommand{\ord}{\operatorname{Ord}}

\newcommand{\lesdot}{\mathrel{\mathord{<}\!\!\raise 0.8
pt\hbox{$\scriptstyle\circ$}}}  

\newcommand{\conc}{^\frown\!}

\newcount\skewfactor
\def\mathunderaccent#1#2 {\let\theaccent#1\skewfactor#2
\mathpalette\putaccentunder}
\def\putaccentunder#1#2{\oalign{$#1#2$\crcr\hidewidth
\vbox to.2ex{\hbox{$#1\skew\skewfactor\theaccent{}$}\vss}\hidewidth}}
\def\name{\mathunderaccent\tilde-3 }

\begin{document}

\title[Covering a function by two continuous functions]{Covering 
a function on the plane by two continuous functions
on an uncountable square - the consistency}
\thanks{Research supported by The Israel Science Foundation administered by
The Israel Academy of Sciences and Humanities. Publication No. 585}
\author{Mariusz Rabus}
\address{Institute of Mathematics, \\
The Hebrew University of Jerusalem,\\
Jerusalem 91904, Israel}
\author{Saharon Shelah}
\email{rabus@@mathstat.yorku.ca \quad
shelah@@math.huji.ac.il}
\maketitle
\begin{abstract}
It is consistent that for every function 
$f:{\bbr}\times {\bbr}
\rightarrow {\bbr}$ there is an uncountable set $A\subseteq {\bbr}$
and two continuous functions $f_{0},f_{1}:D(A)  \rightarrow {\bbr}$
such that $f(\alpha ,\beta)\in \{f_{0}(\alpha ,\beta),f_{1}(\alpha ,\beta)\}$
for every $(\alpha ,\beta) \in A^{2}$, $\alpha\not =\beta$.
\end{abstract}

\section{Introduction}

Suppose that $X$ is a topological space and $f:X\rightarrow {\bbr}$
is a real-valued function on X. Is there a ``large'' subset of $X$
such that the restriction $f\restriction X$ is continuous?
Obviously, if $A\subseteq X$ is a discrete subspace, then
$f\restriction A$ is continuous. Hence in the case when $\dom(f)={\bbr}$,
we can always find an infinite subset on which $f$ is continuous.
The problem whether there is such ``large'' set has been investigated by 
Abraham, Rubin and Shelah in \cite{ARSh}. They proved that it is consistent 
that every function from ${\bbr}$ to ${\bbr}$ is continuous on 
some uncountable set. Later Shelah \cite{473}
showed that every function may be continuous on a non-meager set.

In this paper we consider functions on the plane, ${\bbr}\times 
{\bbr}$. The reasonable question to ask in this case is:
is there a ``large'' set $A\subseteq {\bbr}$ such that on 
$A\times A$ the function $f$ can be cover by two continuous functions? 
Note that we could not hope for $f$ to be just continuous on $A\times A$, e.g.,
if $g$ is a Sierpinski partition, then for every 
uncountable set $A$, $g$ is not continuous on $A\times A$. 
The main result of this paper is the following theorem. For technical
reasons 
we consider squares  without the diagonal, i.e. for a set $A$ we 
consider $D(A)=\{(x,y):x,y \in A,\; x\not =y\}$.

\begin{thm}
 Assume $2^{\aleph_{l}}=\aleph_{l+1}$ for
$l<4$, and  
$\lozenge_{s}(\aleph_{4},\aleph_{1},\aleph_{0})$, see below. 
Then there is a forcing
notion $P$ which preserves  cardinals and cofinalities and such that 
in $V^{P}$, $\con = \aleph_{4}$ and for every function $f:{\bbr}\times {\bbr}
\rightarrow {\bbr}$ there is an uncountable set $A\subseteq {\bbr}$
and two continuous functions $f_{0},f_{1}:D(A)  \rightarrow {\bbr}$
such that $f(\alpha ,\beta)\in \{f_{0}(\alpha ,\beta),f_{1}(\alpha ,\beta)\}$
for every $(\alpha ,\beta) \in D(A)$.
\end{thm}

The proof breaks down into two parts. In 
Section 2, we prove the consistency
of a guessing principle, diamond for systems. Then, is Section 3, we give 
the proof of the theorem.

\begin{remark}
(1) We can replace $\aleph_{0}$ by any $\mu=\mu^{< \mu}$.

\noindent (2) Our main goal was to prove the consistency of the statement in 
the  theorem with $\con < \aleph_{\omega}$. We get $\con = \aleph_{4}$ 
naturally from the proof, but the values $\aleph_{3}$ or $\aleph_{2}$
may be possible.
\end{remark}
%To replace $2$-place functions by $n$-place, for $\mu=\aleph_{0}$,
%we have to follow \cite{288}, we will deal with it elsewhere.

\subsection{Notation}
We use standard set-theoretic notation. Below we list some frequently used 
symbols.
\begin{itemize}
\item For $A,B$ subsets of ordinals of the same order type, 
$OP_{B,A}$ is the order preserving isomorphism from $A$ to $B$.
\item If $C$ is a set of ordinals, then $(C)'$ denotes the set
of accumulation points.
\item Let $\lambda,\chi$ be cardinals, $\chi$ regular. $S^{\lambda}_{\chi}=
\{\alpha \in\lambda: \cf(\alpha)=\chi\}$.
\item For a statement $\phi$ we define $TV(\phi)=0$ if $\phi$ is true,
otherwise $TV(\phi)=1$.
\item ${\bbr}=\od$.
\item If $M$ is a model, $X \subseteq M$, then $Sk(X)$ is the 
Skolem hull of $X$ in $M$.
\item ${\cal L}[\kappa,\theta)$ is a `universal' vocabulary of cardinality 
$\kappa^{< \theta}$, arity $< \theta$. 
\end{itemize}

\section{Diamond for Systems}

In this section we prove the consistency of 
a guessing principle, diamond for systems 
$\lozenge_{s}$. 

\begin{definition} A sequence 
$\bar{M}=\langle M_{u}:u\in [B]^{\leq 2} \rangle$ is a system of models (of
some fixed language) if:
\begin{enumerate}
\item[(1)] $M_{u} \subseteq \ord$, $B\subseteq \ord$,
\item[(2)] $B \cap M_{u} =u$ for every $u \in [B]^{\leq 2}$,
\item[(3)] for every $u, v\in [B]^{\leq 2}$, $|u|=|v|$, the models
$M_{u}$ and $M_{v}$ are isomorphic and $OP_{M_{u},M_{v}}$ is the isomorphism
from $M_{v}$ onto $M_{u}$, $OP_{M_{u},M_{v}}(v)=u$,
\item[(4)] for every $u,v \in [B]^{\leq 2 }$, 
$M_{u} \cap M_{v}\subseteq M_{u \cap v}$,
\item[(5)] if $|u|=|v|$, 
$u'\subseteq u$, $v' =\{\alpha \in v:(\exists \beta \in u')
(|\beta \cap u|=|\alpha \cap v|)\}$, then $OP_{M_{u'},M_{v'}}
\subseteq OP_{M_{u},M_{v}}$, and $OP_{M_{u},M_{u}}=id_{M_{u}}$, 
and if $|w|=|u|$, then
$OP_{M_{u},M_{v}}\circ OP_{M_{v},M_{w}} =OP_{M_{u},M_{w}}$.
\end{enumerate}
\end{definition}

\begin{remark}
See  \cite{289} on the existence of ``nice'' systems of models
for 
$\lambda$ a sufficiently large cardinal, e.g., measurable. Here we 
do not use large cardinals, and try to get a model in which the 
continuum is small, i.e., less than $\aleph_{\omega}$. For this we need 
a suitable guessing principle.
\end{remark}

\begin{definition}[Diamond for systems
$\lozenge_{s}(\lambda,\sigma,\kappa,\theta)$] Let 
$\{C_{\alpha}:\alpha \in \lambda\}$ be a square sequence on $\lambda$.
$\langle \bar{M}^{\alpha}: \alpha \in W \rangle$ is a $\lozenge_{s}(\lambda,
\sigma, \kappa,\theta)$ sequence, (or
$\lozenge_{s}(\lambda,\sigma,
\kappa,\theta)$-diamond for systems) if:
\begin{enumerate}
\item[(A)] $W\subseteq \lambda$ and for every $\alpha \in W$, 
$\bar{M}^{\alpha}=\langle M^{\alpha}_{u}:u\in [B_{\alpha}]^{\leq 2}
\rangle$ is a system of models, 
 $M^{\alpha}_{u}$ is a model of cardinality $\kappa$,
universe $\subseteq \alpha$, vocabulary of cardinality $\leq \kappa$,
arity $< \theta$, a subset of ${\cal L}[\kappa, \theta)$.
\item[(B)] $B_{\alpha} \subseteq \alpha =\sup(B_{\alpha})$,
$\otp(B_{\alpha})=\sigma$, so $\sigma=\cf(\alpha)$. 
\item[(C)] if $M$ is a model with universe $\lambda$, vocabulary of 
cardinality $\leq \kappa$, arity $< \theta$, a subset of 
${\cal L}[\kappa,\theta)$, then for stationarily many 
$\alpha \in W$ for all $u \in [B_{\alpha}]^{\leq 2}$, $M^{\alpha}_{u} \prec
M$,
\item[(D)] if $\alpha , \beta \in W$ and $\otp(C_{\alpha})<
\otp(C_{\beta})$, then 
\begin{enumerate}
\item[(i)] for some $\zeta \in B_{\beta}$, 
$\bigcup\{M^{\beta}_{u}:u\in [B_{\beta}]^{\leq 2} \}-\bigcup\{M^{\beta}_{u}:
u \in [B_{\beta} \cap \zeta]^{\leq 2}\}$ is disjoint from
$\bigcup\{M^{\alpha}_{u}:u \in [B_{\alpha}]^{\leq 2}\}$,
\end{enumerate}
\item[(E)] if $\alpha \not = \beta$ in $W$, $\otp(C_{\alpha})=
\otp(C_{\beta})$, then there is a one-to-one map $h$ from 
$\bigcup_{u \in [B_{\alpha}]^{\leq 2}}M^{\alpha}_{u}$ onto 
$\bigcup_{u \in [B_{\beta}]^{\leq 2}}M^{\beta}_{u}$, order preserving,
mapping $B_{\alpha}$ onto $B_{\beta}$, $M^{\alpha}_{u}$ onto 
$M^{\beta}_{h(u)}$ which is the identity on the intersection of these sets and
the intersection is an initial segment of 
$\bigcup_{u \in [B_{\alpha}]^{\leq 2}}M^{\alpha}_{u}$ and
$\bigcup_{u \in [B_{\beta}]^{\leq 2}}M^{\beta}_{u}$.
\item[(F)] if $\sigma=\kappa$ we may omit $\sigma$.
\end{enumerate}
\end{definition}

\begin{lemma}
Assume: $\kappa < \mu < \lambda$ are uncountable cardinals,
$\lambda = \chi^{+}$, 
$2^{\mu}=\chi$,
$\square_{\lambda}$, $\lozenge_{S^{\chi}_{\sigma}}$,
$\kappa = \kappa^{<\theta}$,
$\mu^{\kappa}=\mu$, $\sigma$, $\chi$, $\kappa$
regular cardinals.

Then there exists a  diamond for systems on $\lambda$, 
$\lozenge_{s}(\lambda,\sigma,\kappa,\theta)$.
\end{lemma}

\Proof Let 
$\bar{C}=\langle C_{\gamma}: \gamma \in \lambda \rangle$ be a square
sequence on $\lambda$. We assume that each $C_{\gamma}$ is closed 
unbounded in $\gamma$,
if $\gamma$ is a limit. 
Let $C_{\gamma}=
\{\alpha_{\zeta}^{\gamma}:\zeta < \otp(C_{\gamma}) \}$. 
First choose a sequence $\langle b^{\alpha}_{i} : i < \chi \rangle$ for 
every $\alpha < \lambda$ such that $b^{\alpha}_{i} \subseteq \alpha$,
$|b^{\alpha}_{i}|< \chi$, $b^{\alpha}_{i}$ increasing, continuous in $i$,
$\alpha = \bigcup \{b^{\alpha}_{i}: i < \chi \}$. Next choose $a_{\alpha}$
for $\alpha < \lambda$ such that 
\begin{enumerate}
\item[(1)] $a_{\alpha} \subseteq \alpha$, 
\item[(2)] if $\cf(\alpha)<\chi$, then $|a_{\alpha}| < \chi$,
\item[(3)] if $\beta \in (C_{\alpha})'$, then $a_{\beta} \subseteq a_{\alpha}$,
\item[(4)] if $\beta  \in C_{\alpha}$ and $i=\otp(C_{\alpha})$, then
$b^{\beta}_{i} \subseteq a_{\alpha}$,
\item[(5)]  if $\otp(C_{\alpha})$ is a limit of limit ordinals, then
$a_{\alpha}= \bigcup_{\beta \in (C_{\alpha})'} a_{\beta}$.
\end{enumerate}

Note that if $\alpha \in S^{\lambda}_{\chi}$, then there is 
a club $C_{\alpha}^{'}\subseteq C_{\alpha}$ such that $\langle a_{\beta}:
\beta\in C_{\alpha}^{'} \rangle$ is an increasing, continuous sequence of 
subsets of $\alpha$ of cardinality $< \chi$ with union $\alpha$.
Let $H_{0},H_{1}$ be functions which witness that $\lambda=\chi^{+}$, i.e.,
$H_{0},H_{1}$ are two place functions, for every $\alpha \in [\chi,\lambda)$,
$H_{0}(\alpha,-)$ is a one-to-one functions from $\alpha$ onto $\chi$
and $H_{1}(\alpha,H_{0}(\alpha,i))=i$ for every $\alpha \in [\chi,\lambda)$
and $i < \alpha$.

Now by induction on $\alpha < \lambda$ we define the truth value of 
`$\alpha \in W$', and if we declare it to be 
true, then  
we also define $\bar{M}^{\alpha}$. Suppose we have defined $W \cap \alpha$
and $\bar{M}^{\beta}$ for $\beta \in W \cap \alpha$. Now consider the 
following properties of an ordinal $\alpha\in \lambda$.
\begin{enumerate}
\item[(a)] $a_{\alpha} \cap \chi =\otp(C_{\alpha})$,
\item[(b)] $a_{\alpha}$ is closed under $H_{0}$ and $H_{1}$,
\item[(c)] for every $\gamma \in a_{\alpha}$ we have:
\begin{enumerate}
\item[(i)] if $\cf(\gamma) < \chi$, then $a_{\alpha} \cap \gamma =
b^{\gamma}_{\otp(C_{\alpha})}$ and $C_{\gamma} \subseteq a_{\alpha}$
and $\otp(C_{\gamma}) \leq \otp(C_{\alpha})$,
\item[(ii)] if $\cf(\gamma)=\chi$, then $\sup(a_{\alpha} \cap \gamma)=
\alpha^{\gamma}_{\otp(C_{\alpha})}$ and $C_{\alpha^{\gamma}_{\otp(C_{\alpha})}}
\subseteq a_{\alpha}$,
\end{enumerate}
\item[(d)] $\cf(\alpha)=\sigma$.
\end{enumerate}
If $\alpha$ does not satisfy one of the 
conditions $(a)$, $(b)$, $(c)$, and $(d)$,
then we declare that $\alpha \not \in W$. So suppose that $\alpha$
satisfies
 $(a)$, $(b)$, $(c)$, and $(d)$. 
Let $\langle M_{\zeta}:\zeta \in\chi\rangle$ 
be the  diamond sequence for $S^{\chi}_{\sigma}$, i.e., 
each $M_{\zeta}$ is a model on $\zeta$,
vocabulary as above, and for every model $M$ on $\chi$,
there are stationarily many $\zeta \in S^{\chi}_{\sigma}$, 
such that $M \cap \zeta
=M_{\zeta}$.  We say that $M_{\zeta}$ is suitable if it is of the form 
$(\zeta , <^{*}_{\zeta} , M^{*}_{\zeta})$, where $<^{*}_{\zeta}$ is a 
well-ordering of $\zeta$.
For each $\zeta$ such that 
$M_{\zeta}$ is suitable, 
let $\xi_{\zeta}= \otp(\zeta , <^{*}_{\zeta})$.
Let $h_{\zeta}:\zeta \rightarrow \xi_{\zeta}$ be the isomorphism
between $(\zeta, <^{*}_{\zeta})$ and $(\xi_{\zeta},<)$.
Let $M^{\oplus}_{\zeta}$ be the model with universe $\xi_{\zeta}$, such that 
$h_{\zeta}$ is the isomorphism between $M^{*}_{\zeta}$ and 
$M^{\oplus}_{\zeta}$. 
For $\alpha\in \lambda$ let $\zeta(\alpha)=\otp(C_{\alpha})$. 
Consider the following properties of $\alpha\in \lambda$.
\begin{enumerate}
\item[(e)] there is a   system $\bar{N}^{\zeta(\alpha)} =\langle 
N^{\zeta(\alpha)}_{s}:
s \in  [\bar{B}_{\zeta(\alpha)}]^{\leq 2} \rangle$, 
$N^{\zeta(\alpha)}_{s} \prec M^{\oplus}_{\zeta(\alpha)}$, 
$||N^{\zeta(\alpha)}_{s}||=\kappa$, 
$\bar{B}_{\zeta(\alpha)}$ cofinal in $\xi_{\zeta(\alpha)}$, 
$\otp(\bar{B}_{\zeta(\alpha)})=\sigma$,
\item[(f)] $\otp(a_{\alpha})=\xi_{\zeta(\alpha)}$.
\end{enumerate}
If $\alpha$ does not satisfy $(e)$, and $(f)$, then declare 
$\alpha \not \in W$. So assume that $\alpha$ satisfies $(e)$ and $(f)$.
Let $g_{\alpha}:\xi_{\zeta(\alpha)} \rightarrow a_{\alpha}$ be the 
order preserving isomorphism. Let ${\bar{M}}^{\alpha}=
\langle M^{\alpha}_{u}:u\in [B_{\alpha}]^{\leq 2}
\rangle$ be the system
of models on $a_{\alpha}$, which is isomorphic to  
$\bar{N}^{\zeta(\alpha)}$ and the isomorphism is $g_{\alpha}$.
If this system satisfies:
\begin{enumerate}
\item[(g)] for every $\beta\in (C_{\alpha})'$ there is $\nu\in B_{\alpha}$
such that $a_{\beta} \cap \bigcup \{M^{\alpha}_{u}:u\in [B_{\alpha}]^{\leq 2}
\} \subseteq \bigcup \{M^{\alpha}_{u}:u\in [B_{\alpha}\cap \nu]^{\leq 2}
\}$,
\end{enumerate}
then we declare $\alpha \in W$. 
This finishes the definition of the diamond for systems sequence,
$\langle \bar{M}^{\alpha}:
\alpha \in W \rangle$.

We have to prove that it is as required. Clauses (A) and (B) are clear.

\medskip
\noindent {\it Proof of clause (C).} 
We need  the following fact, it is proved essentially in 

\noindent \cite{300F},
but for completeness we give the proof at the end of the section.

\begin{lemma} \label{l}
Assume:
\begin{enumerate}
\item[(1)] $\lambda 
=(2^{\mu})^{+}$, $\mu =\mu^{\kappa}$, $\kappa =\cf(\kappa)
>\aleph_{0}$, $\kappa^{<\theta}=\kappa$,
\item[(2)] $M$ is a model with universe $\lambda$, at most $\kappa$ functions
each with $<\theta$ places and $\leq \kappa$ relations including the 
well-ordering of $\lambda$.
\end{enumerate}
Then: for some club $E$ of $\lambda$ for every $\delta \in E$ of cofinality
$\geq \mu^{+}$ we can find $I\subseteq \delta =\sup(I)$ 
and $\langle N_{t}:t\in [I]^{\leq 2},s \in I \rangle$
such that:
\begin{enumerate}
\item[($\alpha$)] $\langle N_{t}:t\in [I]^{\leq 2} \rangle$ is a system of
elementary submodels of $M$, $||N_{t}||=\kappa$.
\end{enumerate}
\end{lemma}

Suppose that $\cal A$ is a model on $\lambda$, $C$ a club 
on $\lambda$. We have to find $\alpha \in C\cap W$ such that 
 $M^{\alpha}_{u} \prec {\cal A}$ for every $u\in [B_{\alpha}]^{\leq 2}$.
Let  $E \subseteq \lambda$ be the
 club given by Lemma \ref{l}. W.l.o.g. we can assume that 
$E \subseteq C'$, where $C'$ is the set of limit points of $C$, (so 
if $\delta \in E$, then $C\cap \delta$ is a club in $\delta$).
Fix $\delta \in S^{\lambda}_{\chi} \cap E$.  Let $f_{\delta}:\delta 
\rightarrow \chi$ be a bijection and let
\begin{equation*} D_{1}=\{\zeta < \chi : \zeta 
\mbox{ \rm is a limit, } f_{\delta} \mbox{ \rm maps }
a_{\alpha^{\delta}_{\zeta}}
\mbox{ \rm onto }\zeta \}.
\end{equation*}
 $D_{1}$ is a $\sigma$-club, i.e., unbounded, closed under 
$\sigma$-sequences.
Let ${\cal A}^{[\delta]}$ be $(\chi , f_{\delta}^{''}(< \restriction \delta),
f_{\delta}^{''}({\cal A}\restriction \delta))$. 
Note that by Lemma \ref{l} we have a system of submodels on ${\cal A}
\restriction \delta$, we transfer this system on ${\cal A}^{[\delta]}$ by
the bijection $f_{\delta}$ and, choosing a subsystem if necessary, we 
can assume that we have an end-extension system on ${\cal A}^{[\delta]}$
which is cofinal in $\chi$, i.e.,
we have $\bar{N}^{*}=\langle N^{*}_{u}:u\in I \rangle$, $I\subseteq \chi$, 
$\sup(I)=\chi$, $N^{*}_{u} \prec {\cal A}^{[\delta]}$ and 
if $\xi < \zeta$ in $I$, then $\min(N^{*}_{\{\zeta\}}\setminus
N^{*}_{\emptyset})>
\sup(N^{*}_{\{\xi\}})$, and if 
$u$ is an initial segment of $v$, then $N^{*}_{u}$ is an initial
segment of $N^{*}_{v}$.
 Hence the set
\begin{equation*}
D_{2}=\{\zeta<\chi: \bigcup_{u \in [\zeta\cap I]^{\leq 2}}N^{*}_{u}
\subseteq  \zeta\}
\end{equation*}
is a club of $\chi$ and such that for every $\zeta\in D_{2}$ there
is a system of models on $\zeta$, ($\langle N^{*}_{u}:u\in [\zeta\cap I]
^{\leq 2} \rangle$).

\noindent Note that the set 
\begin{equation*}
D_{3}=\{\zeta < \chi:\alpha^{\delta}_{\zeta}\in C \text{ and } 
\alpha^{\delta}_{\zeta} 
\mbox{ \rm satisfies conditions }(a)-(d)\}
\end{equation*}
is a $\sigma$-club of $\chi$.
Note that  
${\cal A}^{[\delta]}$ is a model on $\chi$. 
Hence by $\lozenge_{S^{\chi}_{\sigma}}$, 
for stationary many $\zeta \in S^{\chi}_{\sigma}$ we 
have guessed it, i.e., the set \[S=
\{\zeta \in S^{\chi}_{\sigma}: M_{\zeta}={\cal A}^{[\delta]}
\restriction \zeta \}\] is stationary.
Now if
$\zeta \in S\cap (D_{1})' \cap D_{2} \cap D_{3}$  
then $\alpha^{\delta}_{\zeta}\in C$,
and  $\alpha^{\delta}_{\zeta}$
satisfies conditions $(a)$-$(d)$.
Note that $\zeta(\alpha^{\delta}_{\zeta})=\otp(C_{\alpha^{\delta}_{\zeta}})=
\zeta$. Moreover, as $\zeta\in D_{1}\cap S$ we have $\xi_{\zeta}=\otp(a_{
\alpha^{\delta}_{\zeta}})$, i.e., condition $(f)$ holds.
By the construction it follows that condition $(e)$ holds, (the system 
of submodels on $\xi_{\zeta}$ is isomorphic to the system on $a_{
\alpha^{\delta}_{\zeta}}$ given by Lemma \ref{l}). 
Finally, $(g)$ holds, as $\zeta\in (D_{1})'$ and the system of models
of ${\cal A}^{[\delta]}$ is end-extending.

Hence $\alpha^{\delta}_{\zeta}\in W \cap C$, 
and ${\bar{M}}^{\alpha^{\delta}_{\zeta}}$
is a system of models as required.

\medskip
\noindent{\it Proof of clause (E).}
Suppose $\alpha ,\beta \in W$, $\xi=\otp(C_{\alpha})=\otp(C_{\beta})$.
By the construction, both $a_{\alpha}$ and $a_{\beta}$ are isomorphic
to $M^{\oplus}_{\xi}$ and the  isomorphisms are order preserving functions.
Hence $a_{\alpha}$ is order isomorphic to $a_{\beta}$.
Note that $a_{\alpha}\cap \chi =a_{\beta}\cap \chi=\xi$. Since 
both $a_{\alpha}$ and $a_{\beta}$ are closed under $H_{0}$ and $H_{1}$
it follows that $a_{\alpha}\cap a_{\beta}$ is an initial
segment of both $a_{\alpha}$ and $a_{\beta}$. 

\medskip
\noindent{\it Proof of clause (D).}
Suppose that $\alpha , \beta \in W$ and $\otp(C_{\alpha}) <
\otp(C_{\beta})$.
As above, since $a_{\alpha}$ and $a_{\beta}$ are closed under 
$H_{0}$ and $H_{1}$, it follows that $a_{\alpha} \cap a_{\beta}$ is 
an initial
segment of $a_{\alpha}$.  
Let $\gamma = \sup(a_{\alpha} \cap a_{\beta})$. We have four cases, we will 
show that the first three never occur.

\noindent{\it Case 1.} $\gamma \in a_{\alpha} \cap a_{\beta}$.
We can assume that each $a_{\alpha}$ is closed under successor, so this 
case can never happen.

\smallskip
\noindent{\it Case 2.} $\gamma \in a_{\alpha}- a_{\beta}$.
Note that 
$C_{\gamma} \subseteq a_{\alpha}$. Let $\gamma_{1}=\min(a_{\beta}-
\gamma)$. By (c)(i) for $a_{\beta}$ 
it follows that we must have $\cf(\gamma_{1})=\chi$.
Now by (c)(ii), $\gamma =\sup(a_{\beta} \cap \gamma_{1})=
\alpha^{\gamma_{1}}_{\otp(C_{\beta})}$. So $\gamma \in C_{\gamma_{1}}$ and
$\otp(C_{\gamma})=\otp(C_{\beta})$.
Note that $\cf(\gamma)<\chi$. Hence by (c)(i) for $a_{\alpha}$ we
have $\otp(C_{\gamma})\leq \otp(C_{\alpha})$, a contradiction.

\smallskip
\noindent{\it Case 3.} $\gamma \not \in (a_{\alpha} \cup a_{\beta})$.
Let $\gamma_{0}=\min(a_{\alpha}-\gamma)$ and ,$\gamma_{1}
=\min(a_{\beta}-\gamma)$. As above we have $\otp(C_{\gamma})=
\otp(C_{\alpha})$ and $\otp(C_{\gamma})= \otp(C_{\beta})$, a contradiction.

\smallskip
\noindent{\it Case 4.}  
$\gamma \in a_{\beta}- a_{\alpha}$.
 Let $\gamma_{0}=\min(a_{\alpha}-
\alpha)$. We have $\cf(\gamma_{0})=\chi$ and 
$\otp(C_{\gamma})=\otp(C_{\alpha})$, so $C_{\gamma} \subseteq a_{\alpha}$.
Note that $a_{\alpha} \cap \gamma =\bigcup_{\zeta \in C_{\gamma}}(a_{\alpha}
\cap \zeta)$. But for $\zeta\in a_{\alpha}$ with $\cf(\zeta)
<\chi$ we have $a_{\alpha}
\cap \zeta=b^{\zeta}_{\otp(C_{\alpha})}$. Hence
 $a_{\alpha} \cap \gamma =\bigcup_{\zeta \in (C_{\gamma})'}
b^{\zeta}_{\otp(C_{\alpha})} \subseteq a_{\beta_{1}}$, for some 
$\beta_{1}\in (C_{\beta})'$ large enough. Hence by $(g)$ in the definition of 
the diamond for systems sequence, the conclusion follows.

\medskip
\noindent {\bf Proof of Lemma \ref{l}} We prove slightly more. In addition
to the sequence  $\langle N_{t}:t\in [I]^{\leq 2} \rangle$ there is 
a sequence $\langle N'_{\{\alpha\}}:\alpha \in I \rangle$ such that:

\begin{enumerate}
\item[($\beta$)] $N_{\{\alpha\}},N'_{\{\alpha\}}$ 
realize the same $L_{\theta ,\theta}$-type 
over $M$, for $\alpha \in  I$,
\item[($\gamma$)] we have $N'_{\{\alpha\}} \prec N_{\{\alpha\}}$ for 
$\alpha \in I$ and 
for $\alpha < \beta$ in $I$ we have $N_{\{\alpha ,\beta\}}=Sk(N_{\{\alpha\}} 
\cup 
N'_{\{\beta\}})$,
\end{enumerate}

\begin{remark}
(1) Note that for $\alpha < \beta$, $N_{\{\beta\}}$ is not necessarily
a subset of $N_{\{\alpha,\beta\}}$.

\noindent (2) The idea of the proof is to define $N^{*}_{\{0\}}$,
${N^{*}_{\{1\}}}'$ and $N^{*}_{\{0,1\}}$ (and more, see definition of a witness
below).
Then we use it as a blueprint and ``copy'' it many times using 
elementarity, to obtain 
a suitable  system. 
\end{remark}

We can assume that  $M$ has Skolem functions, even for 
$L_{\theta ,\theta}$.
Let $\chi^{*}$ be large enough. Let for $i < \lambda$, ${\cal B}_{i} \prec
(H(\chi^{*}),\in ,<^{*}_{\chi^{*}})$ such that $||{\cal B}_{i}||=2^{\mu}
<\lambda$, and $M \in {\cal B}_{i}$, $ {\cal B}_{i}$ increasing
continuous with $i$, and if $\cf(i) \geq \mu^{+}$ or $i$ non-limit, then 
${\cal B}_{i} \prec_{L_{\mu^{+},\mu^{+}}} 
(H(\chi^{*}),\in ,<^{*}_{\chi^{*}})$.
Let $E=\{\delta <\lambda :
\delta \mbox{ \rm is a limit and } {\cal B}_{\delta} \cap \lambda =
\delta \}$, it is a club of $\lambda$. 
Fix $\delta\in E \cap S^{\lambda}_{\geq \mu^{+}}$. Note that 
${\cal B}_{\delta}
\prec_{L_{\mu^{+},\mu^{+}}}(H(\chi^{*}),\in,<^{*}_{\chi^{*}})$.

We say that $(N^{*}_{\emptyset},
N^{*}_{\{0\}}, {N^{*}_{\{1\}}}', N^{*}_{\{0,1\}},
\alpha_{0},\alpha_{1})$ is a witness if:
\begin{enumerate}
\item[(1)] $N_{u}^{*} \prec M$, $|N_{u}^{*}|=\kappa$, $N^{*}_{\{0\}}
\cap {N^{*}_{\{1\}}}'=N^{*}_{\emptyset}$, $N^{*}_{\emptyset}, 
N^{*}_{\{0\}} \prec M\restriction {\cal B}_{\delta}$,
$N^{*}_{\{0,1\}}=Sk({N^{*}_{\{1\}}}' \cup N^{*}_{\{0\}})$,
\item[(2)] $N^{*}_{\{1\}} \cap {\cal B}_{\delta}=N^{*}_{\emptyset}$,
$\alpha_{0}\in N^{*}_{\{0\}}-N^{*}_{\emptyset}$,
$\alpha_{1}\in {N^{*}_{\{1\}}}'-N^{*}_{\emptyset}$,
\item[(3)] if $\alpha\in N^{*}_{\{0,1\}} \setminus {N^{*}_{\{1\}}}'$, 
$\beta=\min({N^{*}_{\{1\}}}'\setminus \alpha)$, then $\cf(\beta)\geq \mu^{+}$,
\item[(4)] for every $A\subseteq {\cal B}_{\delta}$, $|A| \leq \mu$
there are $N'_{\{1\}} \prec N_{\{1\}}$ and 
$ N_{\{0,1\}}$ such that 
\begin{enumerate}
\item[(a)] $N'_{\{1\}}, N_{\{0,1\}} \prec M \cap {\cal B}_{\delta}$,
\item[(b)] $N'_{\{1\}}$ is order isomorphic to ${N^{*}_{\{1\}}}'$,
\item[(c)] $N_{\{1\}}$ is order isomorphic to $N^{*}_{\{0\}}$,
\item[(d)] $OP_{N_{\{0,1\}},N^{*}
_{\{0,1\}}}$ is an isomorphism from
$N^{*}_{\{0,1\}}$ onto $N_{\{0,1\}}$ which is the identity on 
${N^{*}_{\{1\}}}'$,
maps $N^{*}_{\{0\}}$ onto $N_{\{0\}}$,
\item[(e)] for $\alpha  \in N^{*}_{\{0,1\}} \setminus {N^{*}_{\{1\}}}'$, 
if $\beta
=\min(N_{\{1\}}'-\alpha)$, then $OP_{N_{\{0,1\}},N^{*}_{\{0,1\}}}(\alpha)
\in \sup(A \cap \beta ,\beta)$,
\end{enumerate}
\end{enumerate}

\begin{claim} 
There is a witness.
\end{claim}

 We can find ${\cal C} \prec_{L_{\mu},L_{\mu} }
(H(\chi^{*}),\in ,<^{*}_{\chi^{*}})$ such that $||{\cal C}||=\mu$,
$^{\kappa}{\cal C} \subseteq {\cal C}$, $\mu +1 \subseteq {\cal C}$ and 
$(M,{\cal B}_{\delta},\delta)\in {\cal C}$.
As ${\cal B}_{\delta} \prec_{L_{\mu^{+},\mu^{+}}} (H(\chi^{*}),\in ,
<^{*}_{\chi^{*}})$ it follows that there is a function $f$, 
$\dom(f)={\cal C}$, $\rng(f)\subseteq {\cal B}_{\delta}$, $f\restriction 
{\cal C} \cap {\cal B}_{\delta}$ is the identity, $f$ preserves satisfaction 
of $L_{\mu^{+},\mu^{+}}$ formulas, i.e. $f$ is an isomorphism.

Let ${\cal N} \prec (H(\chi^{*}),\in ,<^{*}_{\chi^{*}})$ be such that 
$\{{\cal B}_{\delta},{\cal C},f,\delta\}\in {\cal N}$, $||{\cal N}||=
\kappa$. Let ${\cal N}_{1}={\cal N} \cap {\cal C}$, ${\cal N}_{0}=
{\cal N} \cap {\cal B}_{\delta}$. Let  ${\cal N}_{0}'=f({\cal N}_{1})$,
note that ${\cal N}_{0}' 
\subseteq {\cal N}_{0}$. Let $\delta_{0}=
f(\delta_{1})$.
W.l.o.g. we can assume that ${\cal N}=
Sk({\cal N}_{0},{\cal N}_{1})$. Let ${\cal N}_{\emptyset}=
{\cal B}_{\delta}\cap {\cal C}\cap {\cal N}$. 
We claim that $({\cal N}_{\emptyset},
{\cal N}_{0},{\cal N}_{1}',{\cal N},
\delta_{0},\delta_{1})$
is a witness.
Note that 
\begin{enumerate}
\item[$(*)$] if $\alpha \in {\cal N} \cap (\delta +1)$, then
$\min({\cal C}-\alpha) \in {\cal N}_{1}$.
\end{enumerate}
Let us check condition (3). Suppose that $\alpha \in {\cal N}-{\cal N}_{1}$
and let $\beta =\min({\cal N}_{1}-\alpha)$. Note that by $(*)$ we have
$\beta =\min({\cal C}-\alpha)$. But as $\mu +1 \subseteq {\cal C}$
and ${\cal C} \prec (H(\chi^{*}),\in ,<^{*}_{\chi^{*}})$ we must have
$\cf(\beta) \geq \mu^{+}$.

Now to verify (4), suppose that there is a set $A$ such that the conclusion 
of (4) fails. Then $A$ is definable from: 
${\cal N}_{1}$, the isomorphism type of ${\cal N}$ over ${\cal N}_{1}$
and the isomorphism type of ${\cal N}_{0}$ over ${\cal N}_{0}'$.
As ${\cal N}_{1}$, ${\cal N}_{\emptyset}$ are in 
${\cal C}$ and ${\cal C} \prec_{L_{\mu},L_{\mu}} 
(H(\chi^{*}),\in ,<^{*}_{\chi^{*}})$ and $\kappa < \mu$ it follows
that such set $A$ is in ${\cal C}$. But now the witness itself is a 
counterexample. Note that clause (e) follows  from $(*)$.

\begin{claim} If there is a witness, 
then there is a system as required,
(for our $\delta \in E \cap S^{\lambda}_{\geq \mu^{+}}$).
\end{claim}

By induction on $\alpha < \mu^{+}$ we define $\delta_{\alpha}<\delta$ and 
a system $\langle
N_{\{\alpha\}}',N_{\{\alpha\}},
N_{\{\alpha ,\beta\}}\rangle$, for $\beta < \alpha$.

Suppose that we have defined the system for all $\beta < \alpha$.
Let $A=\bigcup\{N_{u}:u\in [\{\delta_{\beta}:\beta < \alpha\}]^{\leq 2}\}$.
Let $N_{\{\alpha\}}'$ and $N_{\{\alpha\}}$, $N_{\{0,\alpha\}}$ 
be as in the definition of a witness,
for the above $A$. For $\beta < \alpha$ let $N_{\{\beta ,\alpha\}}=
Sk(N_{\{\beta\}}, N_{\{\alpha\}}')$. It follows that $N_{\alpha}$ is 
isomorphic to ${\cal N}_{0}$ and $N_{\{\beta ,\alpha\}}$ is isomorphic to 
${\cal N}$. 
Let $\delta_{\alpha}= OP_{N_{\{0,\alpha\}},N^{*}_{\{0,1\}}}(\alpha_{0})$.
Note that $I=\{\delta_{\alpha}: \alpha < \mu^{+}\}$ is such that 
$\sup(I)=\delta$ and $N_{u} \cap I=u$ for every $u \in [I]^{\leq 2}$.
This finishes the proof.

\section{Proof of the Theorem}

Start with a model satisfying the assumptions of the theorem, i.e., we have 
$2^{\aleph_{l}}=\aleph_{l+1}$ for $l<4$, $\{C_{\alpha}:\alpha\in 
\omega_{4}\}$ is a square sequence and $\langle 
\bar{M}^{i}:
i \in W \rangle$ is a diamond for systems,  
$\lozenge_{s}(\aleph_{4},\aleph_{1},\aleph_{1},\aleph_{0})$. Let
$\bar{M}^{i}=\langle M^{i}_{u}:u\in [\bar{B}_{i}]^{\leq 2}\rangle$ and
let $\bar{B}_{i}=\{\alpha^{i}_{\epsilon}: 
\epsilon < \omega_{1}\}$ be the increasing enumeration.

\begin{definition} (1) A set $b \subseteq \alpha$ is 
$\bar{Q}\restriction \alpha$-closed, i.e. $\alpha \in b
\Rightarrow a_{\alpha} \subseteq b$.

\noindent (2) 
$\cal{K}={\cal K}_{\mu}$ is the 
family of FS-iterations 
$\bar{Q}=
\langle P_{\alpha}, Q_{\alpha}, a_{\alpha},: \alpha < \alpha^{*}
\rangle$ such that:
\begin{enumerate}
\item[(a)] $a_{\alpha} \subseteq \alpha$,
\item[(b)] $|a_{\alpha}| \leq \mu$,
\item[(c)] $\beta \in a_{\alpha} \Rightarrow a_{\beta} \subseteq 
a_{\alpha}$,
\item[(d)] for $b \subseteq \alpha$, $P^{*}_{b}=\{p\in P_{\alpha}:
 \dom(p) \subseteq b \mbox{ \rm and }(\forall \beta \in \dom(p))
p(\beta) \mbox{ \rm is a }$ $ P^{*}_{b \cap \alpha} \mbox{ \rm name }\}$,
\item[(e)] $Q_{\alpha}$ is a $P^{*}_{a_{\alpha}}$-name, (see 3.2 below),
\item[(f)]  $P^{*}_{\alpha^{*}}$  has the property K, (= Knaster).
\end{enumerate}
\end{definition}

\begin{remark}
The above definition proceeds by induction on $\alpha^{*}$, so 
part (d) is not circular.
\end{remark}

\begin{lemma}
Suppose 
$\bar{Q}=
\langle P_{\alpha}, Q_{\alpha}, a_{\alpha},: \alpha < \alpha^{*}
\rangle \in {\cal K}$. If $b\subseteq \alpha^{*}$ is $\bar{Q}$-closed,
then $P^{*}_{b} \lesdot P^{*}_{\alpha^{*}}$.
\end{lemma}

\Proof Straightforward, see \cite{288} and \cite{289}.

\medskip

 Let $f: \omood \rightarrow \aleph_{1}$ be 
one-to-one, such that if 
$\eta \vartriangleleft \nu$, then $f(\eta) \vartriangleleft f(\nu)$.
For $\rho \in
 \omod$ let $w_{\rho}=\{f(\rho\restriction i):
i<\aleph_{1}\} 
\in [\aleph_{1}]^{\aleph{1}}$. Note that if $\rho_{1} \not = \rho_{2}$ in 
$\omod$,
then $|w_{\rho_{1}} \cap w_{\rho_{2}}| < \aleph_{1}$. 
Let $R$ be the countable support forcing adding $\aleph_{4}$ many Cohen
subsets of $\omega_{1}$,
$\rho_{i}$
$(i<\omega_{4})$. 
Note that in $V^{R}$, $\{w_{\rho_{i}}:i\in \omega_{4}\}$ is a family of almost
disjoint, uncountable subsets of $\omega_{1}$.
Let $B_{i}=\{\alpha^{i}_{\epsilon}:\epsilon\in w_{\rho_{i}}\}$. Note that 
$\{M^{i}_{u}:u\in [B_{i}]^{\leq 2} \}$ is still a system of models on 
$i$, hence without loss of generality we can assume that
$w_{\rho_{i}}=\omega_{1}$.
For $\zeta\in \omega_{1}$ define $B_{i}(\zeta)=\{
\alpha^{i}_{\epsilon}:\epsilon < \zeta\}$.
In $V^{R}$ we shall 
define an iteration $\langle P_{i}, Q_{i}, a_{i}:i<\chi \rangle \in 
{\cal K}_{\aleph_{4}}$.
Working in $V^{R}$, we define 
$\bar{Q}\restriction i$, 
by induction on $i < \omega_{4}$, 
and we prove that it is as in 3.1 (in $V^{R}$).

We call $i$ good if it satisfies:
 $i \in W$, each $M^{i}_{u}$ has a predetermined predicate describing
$\bar{Q}\restriction M^{i}_{u}$ (as an $R$-name, with the limit
$\name{P}^{i}_{u}$)
 and an $R\restriction M^{i}_{u}*\name{P}^{i}_{u}$-name  
$\name{f}$ for a function from $\od \times 
\od$ into $\od$ and each $M^{i}_{u}$ is 
$\bar{Q}$-closed.  (Recall that we do not distinguish between
the model $M^{i}_{u}$ Nan its universe). In this case we put 
$a_{i}=\bigcup\{M^{i}_{u}:u\in [B_{i}]^{\leq 2}\}$ and define 
$Q_{i}$ below.

If $i$ is not good we put $a_{i}=\emptyset$ and define
 $Q_{i}$ to be the Cohen forcing, i.e., $Q_{i}=(\ood, \vartriangleleft)$.
We can assume that if $\alpha \in B_{i}$, then $Q_{\alpha}$ is Cohen,
(or just replace $B_{i}$ by $\{\alpha+1:\alpha\in B_{i}\}$).
For $\alpha\in B_{i}$, let $r_{\alpha}$ be the Cohen real forced by 
$Q_{\alpha}$.

\begin{remark} The reason we add $\aleph_{4}$ almost disjoint
subsets of $\omega_{1}$ is that, in $V^{R}$, if $i\not =j$ are good and 
$\otp(C_{i})=\otp(C_{j})$,
then the systems associated with $i$ and $j$ are almost 
disjoint, i.e., there is $\zeta\in \omega_{1}$ such that 
\begin{eqnarray*}
\lefteqn{(\bigcup\{M^{i}_{u}:u\in [B_{i}]^{\leq 2}\}) \cap 
(\bigcup\{M^{j}_{u}:u\in [B_{j}]^{\leq 2}\})
\subseteq}  \\
& (\bigcup\{M^{i}_{u}:u\in [B_{i}(\zeta)]^{\leq 2}\}) \cap 
(\bigcup\{M^{j}_{u}:u\in [B_{j}(\zeta)]^{\leq 2}\})
\end{eqnarray*}

Note that if $\otp(C_{i})\not = \otp(C_{j})$ then we have almost 
disjointness by 2.2(D)(i).
\end{remark}

\medskip
\noindent {\bf Notation} For $\xi,\zeta \in \omega_{1}$ let 
$Z^{i}_{\xi,\zeta}=M^{i}_{\{\alpha^{i}_{\xi},\alpha^{i}_{\zeta}\}} 
\cup M^{i}_{\{\alpha^{i}_{\xi}\}}
\cup M^{i}_{\{\alpha^{i}_{\zeta}\}}$, 
$Z^{i}_{\xi}=M^{i}_{\{\alpha^{i}_{\xi}\}}$.

Now we fix a good $i$. Our goal is to define $Q_{i}$.

\begin{definition} For $p,q \in R$ (or in $P_{\omega_{4}}^{*}$), 
$\dom(p),\dom(q) \subseteq  Z^{i}_{0,1}$ 
we say that $p$ and $q$ are dual if
$OP_{Z^{i}_{1},Z^{i}_{0}}(
p\restriction Z^{i}_{0})=q\restriction Z^{i}_{1}$ 
and 
$OP_{Z^{i}_{1},Z^{i}_{0}}(
q\restriction Z^{i}_{0})=p\restriction Z^{i}_{1}$. 
\end{definition}

Using $G_{R\restriction M^{i}_{\emptyset}}$ we choose, by induction
on $k < \omega$, conditions $r^{i}_{\eta}$, $r^{i,l}_{\eta}\in R$ for
$\eta \in \kkd$, $l<2$, such that:
\begin{enumerate}
\item[(a)] $r^{i}_{\eta} \in (R\restriction Z^{i}_{0})/
G_{R\restriction M^{i}_{\emptyset}}$.
\item[(b)] $\nu\vartriangleleft \eta \Rightarrow 
r^{i}_{\nu} \leq r^{i}_{\eta}$.
\item[(c)] if $l = m +1$, if $\eta \in \mmmd$, $l<2$, then 
$r^{i,l}_{\eta}  \in (R\restriction Z^{i}_{0,1})/
G_{R\restriction M^{i}_{\emptyset}}$ and $r^{i}_{\eta} \leq r^{i,l}_{\eta}
\restriction Z^{i}_{0} \leq r^{i}_{\eta \conc <l>}$ and 
$OP_{Z^{i}_{1},Z^{i}_{0}}(r^{i}_{\eta}) \leq
r^{i,l}_{\eta}\restriction Z^{i}_{1} \leq 
OP_{Z^{i}_{1},Z^{i}_{0}}(r^{i}_{\eta \conc <1-l>})$,
and $r^{i,0}_{\eta}$ and $r^{i,1}_{\eta}$ are dual.
\item[(d)] $r^{i,l}_{\eta}$ forces that $A^{\eta,l}_{k}=
\{ \name{p}^{\eta,l}_{k,n}:n\in\omega
 \}$ is a predense subset  of $P^{*}_{Z^{i}_{0,1}}$, 
such that each $\name{p}^{\eta,l}_{k,n}$ forces the
value $f^{\eta,l}_{k,n}$ of $\name{f}(r_{\alpha^{i}_{0}}, r_{\alpha^{i}_{1}})
\restriction k$.
\item[(e)]  
$A^{\eta,0}_{k}$ and $A^{\eta,1}_{k}$ are dual, i.e. for every
$m\in \omega$, $\name{p}^{\eta,0}_{k,m}$ and $\name{p}^{\eta,1}_{k,m}$
are dual. Moreover if $k_{1}<k_{2}$, then 
$A^{\eta,l}_{k_{2}}$ refines $A^{\eta,l}_{k_{1}}$.
\end{enumerate}

Suppose we have $r^{i}_{\eta}$. We define $r^{i,0}_{\eta}$,
$r^{i,1}_{\eta}$ and $A_{k}^{\eta,0}$, $A_{k}^{\eta,0}$ as follows.

$1.$ Let $r_{1}=r^{i}_{\eta} \cup OP_{Z^{i}_{1},Z^{i}_{0}}(r^{i}_{\eta})$.

$2.$ Let $r_{1,0}\geq r_{1}$, $r_{1,0}\in R\restriction Z_{0,1}$,
forces a maximal antichain $A_{1,0}$ of $P^{*}_{Z_{0,1}}$,
such that each element of $A_{1,0}$ forces a  value of 
$\name{f}(r_{\alpha^{i}_{0}},r_{\alpha^{i}_{1}}) \restriction k$.

$3.$ Let $r_{2}= OP_{Z^{i}_{1},Z^{i}_{0}}(r_{1,0}\restriction Z^{i}_{0}) \cup 
 OP_{Z^{i}_{0},Z^{i}_{1}}(r_{1,0}\restriction Z^{i}_{1})$. 
Let $r_{2,1}\geq r_{2}$,
$r_{2,1}\in R\restriction Z_{0,1}$ forces $A_{2,1}$ 
to be a predense subset of $P^{*}_{Z_{0,1}}$ such that each element 
of $A_{2,1}$ forces a value of 
$\name{f}(r_{\alpha^{i}_{0}},r_{\alpha^{i}_{1}}) \restriction k$. Moreover, 
$A_{2,1}=\bigcup\{A_{p}:p\in A_{1,0}\}$, where 
for every $q\in A_{p}$ we have $q \geq 
OP_{Z^{i}_{1},Z^{i}_{0}}(p\restriction Z^{i}_{0}) \cup 
 OP_{Z^{i}_{0},Z^{i}_{1}}(p\restriction Z^{i}_{1})$.

$4.$ Let $r_{3}=OP_{Z^{i}_{1},Z^{i}_{0}}(r_{2,1}\restriction Z^{i}_{0}) \cup 
 OP_{Z^{i}_{0},Z^{i}_{1}}(r_{2,1}\restriction Z^{i}_{1})$.

$5.$ Let $r_{3,0}= r_{3} \cup r_{1,0}$ (note: $r_{3,0}$ is dual to 
$r_{2,1}$). Let $A_{3,0}=\{p \cup OP_{Z^{i}_{1},Z^{i}_{0}}(q\restriction
 Z^{i}_{0}) \cup 
 OP_{Z^{i}_{0},Z^{i}_{1}}(q\restriction Z^{i}_{1}): q \in A_{p}\}$.

$6.$ Let $r^{i,0}_{\eta}=r_{3,0}$, $r^{i,1}_{\eta}=r_{2,1}$,
$A^{\eta,0}_{k}=A_{3,0}$ and $A^{\eta,1}_{k}=A_{2,1}$.

\medskip

Let for $\eta \in \od$, $r^{i}_{\eta}=\bigcup_{k<\omega}
r^{i}_{\eta\restriction k}$.
In $V$  choose $\langle \eta^{*}_{\epsilon} : \epsilon <\omega_{1} \rangle$,
distinct members of $\od$.
Recall that $\rho_{j}$ ($j<\aleph_{4}$) are the Cohen subsets of $\omega_{1}$
forced by $R$
In $V[\langle \rho_{j}:j \in \{i\} \cup a_{i} 
 \rangle]$
we can find $w^{i} \in [\omega_{1}]^{\omega_{1}}$ 
such that 
\begin{enumerate}
\item[($\alpha$)] if $\epsilon \in w^{i}$ then 
$OP_{Z^{i}_{\epsilon
},Z^{i}_{0}}(r_{\eta^{*}_{\epsilon}})
\in G_{R\restriction Z^{i}_{\epsilon}}$,
\item[($\beta$)] if $\epsilon_{0} < \epsilon_{1}$ are in $w^{i}$, 
$l=TV(\eta^{*}_{\epsilon_{0}}<_{lx}
\eta^{*}_{\epsilon_{1}})$, 
then 

$OP_{Z^{i}_{\epsilon_{0},\epsilon_{1}},
Z^{i}_{0,1}}(r^{i,l}_{\eta^{*}_{\epsilon_{0}} 
\cap \eta^{*}_{\epsilon_{1}}}) \in G_{R\restriction 
Z^{i}_{\epsilon_{0}, \epsilon_{1}}}$.
\end{enumerate}
We  choose the members of $w^{i}$ inductively using the fact that 
$R$ has $(<\aleph_{1})$-support.

\medskip
\noindent {\bf Notation} For $\xi \in w^{i}$ denote $r^{i}_{\xi}=
r_{\alpha^{i}_{\xi}}$.

\medskip
\noindent Let $H$ be $R$-generic and $G$ be $P^{*}_{a_{i}}$-generic.
In $V[H][G]$ we  define $Q_{i}$.
A condition in $Q_{i}$ is $(u,v,\bar{\nu},\bar{m},F_{0},F_{1})$, where:
\begin{enumerate} 
\item[(1)] $u$ is a finite subset of $w^{i}$.
\item[(2)] $v$ is a finite set of elements of the form 
$(\eta ,\rho)$, where
\begin{enumerate} 
\item[(a)] $\eta,\rho \in \ood$, $\lh(\eta)=\lh(\rho)$, $\rho \not = \eta$,
\item[(b)] $\eta \vartriangleleft r^{i}_{\alpha}$, 
$\rho \vartriangleleft r^{i}_{\beta}$ for some
$\alpha ,\beta \in u$ and if $\nu=\eta^{*}_{\alpha}\cap\eta^{*}_{\beta}$
then for every $\gamma \in u$ we have: if 
$\eta \vartriangleleft r^{i}_{\gamma}$, then 
$\eta^{*}_{\gamma}\restriction (\lh(\nu)+1) =\eta^{*}_{\alpha}
\restriction (\lh(\nu)+1)$, and if $\rho \vartriangleleft r^{i}_{\gamma}$, then
$\eta^{*}_{\gamma}\restriction (\lh(\nu)+1) =\eta^{*}_{\beta}
\restriction (\lh(\nu)+1)$.
\end{enumerate}
\item[(3)] $\bar{\nu}$ is a function from $v$ into
$\ood$ such that 
for $(\eta ,\rho)\in v$ we have:
$\bar{\nu}(\eta,\rho)$ is such that 
there is  $\alpha,\beta \in u$ such that 
$\eta \vartriangleleft 
r^{i}_{\alpha}$, $\rho \vartriangleleft r^{i}_{\beta}$ and 
$\bar{\nu}(\eta,\rho)=
\eta^{*}_{\alpha}
\cap \eta^{*}_{\beta}$, ($\bar{\nu} $ is well defined by (2)).
\item[(4)] $\bar{m}$ is a function from $v$ to 
$\omega$. For $(\eta,\rho)\in v$, $\bar{m}(\eta,\rho)$ is such that 
for every $\alpha,\beta\in u$ such that 
$\eta \vartriangleleft r^{i}_{\alpha}$, $\rho \vartriangleleft r^{i}_{\beta}$, 
we have $OP_{Z^{i}_{\alpha,\beta},
Z^{i}_{0,1}}(p^{\nu,l}_{lh(\eta),\bar{m}(\eta,\rho)})\in G$, 
where $l=TV(\eta^{*}_{\alpha}<_{lx}
\eta^{*}_{\beta})$ and $\nu =\eta^{*}_{\alpha} \cap\eta^{*}_{\beta}$.
\item[(5)] For $l=0,1$, $F_{l}$ is a function from $v$
into $\ood$, defined by:
for $(\eta,\rho)\in v$, $F_{l}(\eta,\rho)$ is the value
of $\name{f}(r_{0},r_{1})\restriction lh(\eta)$ forced by
$p^{\bar{\nu}(\eta ,\rho),l}_{lh(\eta),\bar{m}(\eta,\rho)}$.
\item[(6)] For $(\eta,\rho), (\eta_{1},\rho_{1})\in v$, if $\eta 
\vartriangleleft \eta_{1}$ and $\rho \vartriangleleft \rho_{1}$,
then $F_{l}(\eta,\rho) \vartriangleleft F_{l}(\eta_{1},\rho_{1})$, for 
$l=0,1$.
\end{enumerate}
Order: $(u,v,\bar{\nu},\bar{m},F_{0},F_{1}) \leq 
(u^{1},v^{1},\bar{\nu}^{1},\bar{m}^{1},F_{0}^{1},F_{1}^{1})$ if 
\begin{enumerate}
\item[(7)] $u\subseteq u^{1}$,
\item[(8)] $v\subseteq v^{1}$,
\item[(9)] $F_{l}=F_{l}^{1} \restriction v$, 
$\bar{\nu}=\bar{\nu}^{1}\restriction v$,
$\bar{m}=\bar{m}^{1}\restriction v$, 
$l=0,1$.
\end{enumerate}

\begin{lemma} Suppose $(q_{\alpha},p_{\alpha})$, (for $\alpha \in \omega_{1}$),
are in $P^{*}_{a_{i}}*Q_{i}$, $q_{\alpha}$ forces $p_{\alpha}$
to be a real $6$-tuple in $Q_{i}$, not just a $P^{*}_{a_{i}}$-name
of such a tuple, $\dom(q_{\alpha})$
$(\alpha \in \omega_{1})$ form a delta system with the 
root $\Delta$, $\zeta\in \omega_{1}$. Let $b=
\bigcup
\{M^{i}_{u}:u \in [B_{i}(\zeta)]^{\leq 2}\} $.
Suppose
$\Delta-\{i\} \subseteq b$ and 
$\dom(q_{\alpha}) \cap b=
\Delta$ for $\alpha \in \omega_{1}$.

Then there is an 
uncountable set $E\subseteq \omega_{1}$ such that for every
$\alpha,\beta\in E$, $(q_{\alpha},p_{\alpha})$ and 
$(q_{\beta},p_{\beta})$ are compatible, 
moreover if $q\in P^{*}_{b}$, $q \geq q_{\alpha}\restriction b,q_{\beta}
\restriction b$, then $q, (q_{\alpha},p_{\alpha})$
and $(q_{\beta},p_{\beta})$ are compatible.
\end{lemma}

\Proof 
By thinning out we can find an uncountable set $E\subseteq \omega_{1}$
such  that:
\begin{enumerate}
\item[(a)] For $\alpha \in E$
let $w_{\alpha}=\bigcup\{u \in [B_{i}]^{\leq 2}:\dom(q_{\alpha})\cap
M^{i}_{u}\not =\emptyset\}$, (each $w_{\alpha}$ is finite).
The sets $w_{\alpha}$, $(\alpha\in E)$ form a delta
system with the  root $w$ and if $\alpha<\beta$, $\xi \in w_{\alpha},
\zeta \in w_{\beta}$, then $\xi\leq \zeta$.
\item[(b)] $u^{p_{\alpha}}$ $(\alpha\in E)$ form
a delta system with the root $u$ and $\alpha<\beta$, $\xi \in u^{p_{\alpha}},
\zeta \in u^{p_{\beta}}$, then $\xi\leq \zeta$,  $|u^{p_{\alpha}}|=n^{*}$.
\item[(c)] $v^{p_{\alpha}}=v^{*}$ for $\alpha\in E$
and the structures $(u^{p_{\alpha}}, \{q_{\alpha}(\xi):
\xi \in u^{p_{\alpha}}\},
v^{*},\{\eta^{*}_{\xi}\restriction m^{*}
:\xi \in u^{p_{\alpha}}\})$ are isomorphic, 
(isomorphism given by the order preserving bijection between respective
$u^{p_{\alpha}}$'s),
where 
$m^{*}$ is such that $\lh(\eta^{*}_{\xi}\cap \eta^{*}_{\zeta}) <m^{*}$
for every $\xi\not =\zeta$ in $u^{p_{\alpha}}$.
\end{enumerate}

\begin{lemma}
$P_{i+1}$ has the property $K$.
\end{lemma}

\Proof 
Let $\{p_{\alpha}:\alpha\in  \omega_{1}\}$ be an uncountable subset
of $P_{i+1}$. W.l.o.g. we can  assume that $\dom(p_{\alpha})$,
$(\alpha\in\omega_{1})$ form a delta system  with  the root  
$\Delta$. We have to find an    uncountable subset $E\subseteq
\omega_{1}$  such that for any $\alpha  ,\beta \in E$, 
$p_{\alpha}$ and $p_{\beta}$ are compatible.
We prove it by induction on $k=|\Delta|$.

For $k=0$, trivial. For  the induction step assume that
$\Delta=\{i_{0},\ldots ,i_{k}\}$ ordered by $\vartriangleleft$,
where for $\alpha ,\beta < \omega_{4}$, we define $\alpha \vartriangleleft
\beta$ iff $\otp(C_{\alpha}) < \otp(C_{\beta})$ 
or $\otp(C_{\alpha})=\otp(C_{\beta})$ and 
$\alpha < \beta$.

By the induction hypothesis there is an uncountable set $E' \subseteq 
\omega_{1}$ such that for $\alpha ,\beta \in E'$, 
$p_{\alpha}\restriction \bigcup_{l<k}a_{i_{l}}$ and 
$p_{\beta}\restriction \bigcup_{l<k}a_{i_{l}}$  are compatible.
Note that 
there is $\zeta\in \omega_{1}$ such that 
$a_{i_{k}} \cap (\bigcup_{l<k}a_{i_{l}}) \subseteq
\bigcup\{M^{i_{k}}_{u}:u\in [B_{i_{k}}(\zeta)]^{\leq 2}\}$,
(see 2.2(D)).
Now use the  previous lemma.

\medskip
Now suppose that $G(i)$ is $Q_{i}$-generic. Let 
\[A'=\bigcup\{u:\exists (v,\bar{\nu},\bar{m},F_{0},F_{1}), 
(u,v,\bar{\nu},\bar{m},F_{0},F_{1})\in G(i)\}.\]
In $V[G]$ let 
$A=\{r^{i}_{\alpha}:\alpha\in A'\}$ and let   
$f_{l}:[A]^{2}\rightarrow \od$ be defined by:
\begin{eqnarray*}
\lefteqn{f_{l}(r^{i}_{\alpha},r^{i}_{\beta})=\bigcup\{F_{l}(\eta,\rho):
\exists (u,v,\bar{\nu},\bar{m},F_{0},F_{1})\in G(i),}  \\
 & \alpha,\beta\in u,
(\eta,\rho)\in v, \eta \vartriangleleft r^{i}_{\alpha},
\rho \vartriangleleft r^{i}_{\beta}\}.
\end{eqnarray*}

Let ${\cal V}
=\bigcup\{v: \exists (u,\bar{\nu},\bar{m},F_{0},F_{1}): 
(u,v,\bar{\nu},\bar{m},F_{0},F_{1})\in G(i)\}$.

\begin{lemma}
\begin{enumerate}
\item[(1)] For every $\alpha,\beta \in A'$ and $n\in \omega$ 
there is $(\eta,\rho)\in {\cal V}$ such that 
$\lh(\eta)=\lh(\rho)\geq n$ and 
$\eta \vartriangleleft r_{\alpha}$ and 
$\rho \vartriangleleft r_{\beta}$,
\item[(2)] $A$ is uncountable,
\item[(3)] $f_{0},f_{1}$ are continuous,
\item[(4)] for every $(\alpha,\beta)\in [A]^{2}$, 
if $l=TV(\eta^{*}_{\alpha}<_{lx}\eta^{*}_{\beta})$, then 
$f(r^{i}_{\alpha},r^{i}_{\beta})=f_{l}(r^{i}_{\alpha},r^{i}_{\beta})$.
\end{enumerate}
\end{lemma}

\Proof (1) and (2) follow by a density argument. 
To prove (1) suppose that $(p,q)\in P_{i}*Q_{i}$,
$p$ forces that $\alpha,\beta \in u^{q}$. W.l.o.g. $\alpha,\beta\in 
\dom(p)$. Let $p_{1}\in P_{i}$
be such that $\dom(p)=\dom(p_{1})$, $p(\zeta)=p_{1}(\zeta)$ for 
$\zeta\in \dom(p)\setminus \{\alpha,\beta\}$, $p(\alpha)
\vartriangleleft p_{1}(\alpha)$,
$p(\beta)
\vartriangleleft p_{1}(\beta)$, $\lh(p_{1}(\alpha))=\lh(p_{1}(\beta))\geq n$,
(remember that $Q_{\alpha}, Q_{\beta}$ are Cohen).
Let $\eta=p_{1}(\alpha)$, $\rho=p_{1}(\beta)$, $\nu=\eta^{*}_{\alpha}
\cap \eta^{*}_{\beta}$, $l=TV(\eta^{*}_{\alpha}<_{lx}\eta^{*}_{\beta})$. 
Let $m\in \omega$ be such that $OP_{Z_{\alpha,\beta},Z_{0,1}}(
p^{\nu,l}_{\lh(\eta),m})$ is compatible with $p_{1}$, and 
let $p_{2}$ be the common upper bound.  Now define $q_{1}\geq q$ as 
follows. $u^{q_{1}}=u^{q}$, $v^{q_{1}}=v^{q}\cup\{(\eta,\rho)\}$,
$\bar{\nu}^{q_{1}}(\eta,\rho)=\nu$, $\bar{m}^{q_{1}}(\eta,\rho)=m$,
$F_{l}^{q_{1}}(\eta,\rho)$ is the value forced by $p^{\nu,l}_{\lh(\eta),m}$.
Hence $(p_{2},q_{1}) \geq (p,q)$ and it forces what is required.

To prove (2) it 
is enough to show, in $V^{R}$, that for every $\alpha\in \omega_{1}$
and $(p,q)\in P_{i}*Q_{i}$
there is $\beta > \alpha$ and $(p_{1},q_{1})\geq (p,q)$, such that
$\beta \in u^{q_{1}}$. Let $\beta > \alpha$
be such that 
$\dom(p) \cap Z^{i}_{\gamma,\beta} \subseteq M^{i}_{\emptyset}$ and $\beta > 
\gamma$ for 
every $\gamma \in u^{q}$. Let $\gamma \in u^{q}$ be such that 
$(\eta^{*}_{\gamma_{1}} \cap \eta^{*}_{\beta}) \vartriangleleft
(\eta^{*}_{\gamma} \cap \eta^{*}_{\beta})$ for every $\gamma_{1}\in u^{q}$.
Define condition $q_{1}(\beta)=q(\gamma)$ and let $p_{1}$ be a condition
extending $p$ and each of conditions $OP_{Z^{i}_{\gamma_{1},\beta},Z^{i}_{0,1}}
(p^{\bar{\nu}(\eta,\rho),l}_{\lh(\eta),\bar{m}(\eta,\rho)})$ such that 
$(\eta,\rho)\in v$, $\eta \vartriangleleft q(\gamma_{1})$,
$\rho \vartriangleleft q(\gamma)$ and 
$l=TV(\eta^{*}_{\gamma_{1}}< \eta^{*}_{\beta})$.
 Finally extend $q$ to $q_{1}$ such that $u^{q_{1}}=u^{q}\cup \{\beta\}$.

Condition (3) 
follows from (1) and (5) and (6) in the definition of $Q_{i}$.

To prove (4) it is enough to 
show that for every $n\in \omega$, 
$f(r^{i}_{\alpha},r^{i}_{\beta})\restriction n
=
f_{l}(r^{i}_{\alpha},r^{i}_{\beta})\restriction n$. 
By condition (1) there is $(\eta,\rho) \in V$ such that $k=\lh(\eta)\geq n$
and $\eta \vartriangleleft r^{i}_{\alpha}$ and 
$\rho \vartriangleleft r^{i}_{\beta}$. 
Recall that  
$p=p^{\bar{\nu}(\eta,\rho),l}_{\lh(\eta),\bar{m}(\eta,\rho)}$ forces
that $f(r^{i}_{0},r^{i}_{1})\restriction k=h$ for some fixed $h$.
Now working in $V$ consider $(r^{i,l}_{\eta^{*}_{\alpha} 
\cap \eta^{*}_{\beta}}, p) \in R*P_{i}\restriction Z^{i}_{0,1}$. By the 
construction the condition $(r',p)=OP_{Z^{i}_{\alpha,\beta},Z^{i}_{0,1}}
(r^{i,l}_{\eta^{*}_{\alpha} \cap \eta^{*}_{\beta}},p)\in H*G$, 
and forces that $f(r^{i}_{\alpha},r^{i}_{\beta})=h$. 
On the other hand, by definition
$F_{l}(\eta,\rho) =h$ and 
$F_{l}(\eta,\rho) \vartriangleleft f_{l}(r^{i}_{\alpha},r^{i}_{\beta})$
This finishes the proof.


\begin{thebibliography}{[Sh 300F]}
\bibitem[ARSh]{ARSh}Uri Abraham, Matatyahu Rubin, and Saharon
Shelah, {\em On the consistency of some partition theorems for 
continuous 
colorings, and the structure of $\aleph_{1}$-dense real order types,}
Ann. Pure. Appl. Logic { \bf 29 } (1985), 123-206.
\bibitem[Sh 473]{473} Saharon Shelah, 
{\em Possibly every real function is continuous on a 
non-meager set,} Publications de L'Institute Mathematique, Nouvelle serie 
tome { \bf 57 } (71), 1995, 47-60.
\bibitem[Sh 300F]{300F}Saharon Shelach, Universal Classes - Classification 
Theory, forthcoming.
\bibitem[Sh 289]{289}Saharon Shelah, {\em Strong Partition Relatons Below
the Power Set: Consistency, Was Sierpinski Right - II?} In:
Proceedings of the conference on Set Theory and its Applications in honor
of A. Hajnal and V.T. Sos, Budapest 1991.
\bibitem[Sh 288]{288}Saharon Shelah, {\em Consistency of positive
partition theorems for graphs and models}, In: Set Theory and 
Its
Applications Conference, York University, J. Steprans and S. Watson
eds, Lecture Notes in Mathematics 1401, Springer-Verlag 1989, 167-193.
\end{thebibliography}
\end{document}